\documentclass[11pt,reqno]{amsart}
\usepackage{curves} 
\usepackage{graphicx} 
\usepackage{amsmath, amsfonts, amssymb, amscd}

\usepackage{amsmath, amsfonts, amssymb, amscd,
 enumerate, enumitem, color, xcolor, cancel,
graphpap,  mathrsfs
}
\let\lpolish\l

\theoremstyle{plain}
\newtheorem{theo}{Theorem}[section]
\newtheorem{lem}[theo]{Lemma}

\newtheorem{cor}[theo]{Corollary}
\theoremstyle{definition}
\newtheorem{rem}[theo]{Remark}

\newtheorem{definition}[theo]{Definition}

\theoremstyle{plain}

\theoremstyle{definition}

\renewcommand{\=}{:=}

\newcommand{\beq}{\begin{equation}}
\newcommand{\eeq}{\end{equation}}
\renewcommand{\a}{\alpha}
\renewcommand{\b}{\beta}

\renewcommand{\d}{\delta}

\newcommand{\g}{\gamma}

\renewcommand{\k}{\kappa}
\renewcommand{\l}{\lambda}

\renewcommand{\t}{\tau}

\newcommand{\z}{\zeta}
\newcommand{\D}{\Delta}

\newcommand{\bB}{\mathbb{B}}
\newcommand{\bC}{\mathbb{C}}

\newcommand{\bR}{\mathbb{R}}

\newcommand{\cA}{\mathcal{A}}

\newcommand{\cC}{\mathcal{C}}

\newcommand{\cL}{\mathcal{L}}

\newcommand{\cZ}{\mathcal{Z}}

\newcommand{\p}{\partial}

\renewcommand{\square}{\kern1pt\vbox
{\hrule height 0.6pt\hbox{\vrule width 0.6pt\hskip 3pt
\vbox{\vskip 6pt}\hskip 3pt\vrule width 0.6pt}\hrule height0.6pt}\kern1pt}

\DeclareMathOperator\Aut{Aut\;}

\DeclareMathOperator\Id{Id}

\newcommand{\wt}{\widetilde}

\newcommand{\be}{\begin{equation}}
\newcommand{\ee}{\end{equation}}

\def\<#1,#2>{\langle\,#1,\,#2\,\rangle}
\newcommand{\arr}{\begin{array}{rlll}}
\newcommand{\ea}{\end{array}}
\newcommand{\bea}{\begin{eqnarray}}
\newcommand{\eea}{\end{eqnarray}}
\newcommand{\bean}{\begin{eqnarray*}}
\newcommand{\eean}{\end{eqnarray*}}

\def\sideremark#1{\ifvmode\leavevmode\fi\vadjust{
\vbox to0pt{\hbox to 0pt{\hskip\hsize\hskip1em
\vbox{\hsize3cm\tiny\raggedright\pretolerance10000
\noindent #1\hfill}\hss}\vbox to8pt{\vfil}\vss}}}

\newcounter{ssig}
\newcounter{ttig}
\title[Regularity of Kobayashi metric]
{Regularity of Kobayashi metric}
\author[G. Patrizio and A.  Spiro]{Giorgio Patrizio
and Andrea Spiro}


\thanks{{\it Acknowledgments}. This research was partially supported by the Project MIUR ``Real and Complex Manifolds: Geometry, Topology and  Harmonic Analysis'' and by GNSAGA of INdAM}

\begin{document}
\dedicatory{Dedicated to Kang-Tae Kim for his sixtieth birthday}
\begin{abstract}  We review some recent   results  on existence and regularity of
  Monge-Amp\`ere exhaustions  on the smoothly bounded strongly pseudoconvex domains,   which admit at least one  such exhaustion of sufficiently high regularity.  A   main consequence of our results  is the fact  that 
 the  Kobayashi  pseudo-metric $\k$ on an appropriare open subset of each of the above  domains  is  actually a  smooth Finsler metric.
 The class of  domains to which our result apply is very large. It includes for instance  all smoothly bounded  strongly pseudoconvex complete circular domains 
  and all their sufficiently small deformations.  \end{abstract}
\maketitle

\null

\section{Introduction}
 In this note, providing the necessary background, we survey some recent results about the existence of regular Monge-Amp\`ere exhaustions which, in turn, imply regularity properties for the Kobayashi metric. More precisely, for domains $D$ admitting a 
smooth  Monge-Amp\`ere exhaustion centered at a point $z_o\in D$ (that is, a  strictly plurisubharmonic  $\cC^0$ exhaustions $\t: \overline D \to [0,1]$,  which is  $\cC^\infty$ at all points, with  only possible 
 exception at the minimum set $\{\t = 0\}=\{z_o\}$,  and such that   $u \= \log \t$  has  a logarithmic singularity at $z_o$ and satisfies  the complex Monge-Amp\`ere equation $(\p \bar \p u)^n = 0$ at all other points), our results show  that there exists an open neighborhood $D' \subset D$ of $z_o$ such that for each  $z \in D'$  there exists an analogous  smooth  Monge-Amp\`ere exhaustion,  centered at such point (\cite{PS3}). One of the main consequence of this result is  that the  Kobayashi  pseudo-metric  on each  such domain is  actually a  smooth complex Finsler metric on an open subset $D' \subset D$.
Since any smoothly bounded  strongly pseudoconvex complete circular domain  
 and any of  its sufficiently small smooth deformations have at least one Monge-Amp\`ere exhaustion and our proof shows who to determine when $D' = D$, our result reveals that there exists a new  large  class of domains, on which 
 the Kobayashi metric has extremely high regularity properties.  In fact, the results in \cite{PS3} are   proven  for closed  strongly pseudoconvex domains with Monge-Amp\`ere exhaustions of  class $\cC^{r, \a}$, $r \geq 4$, $\a > 0$,   and imply  regularity for the Kobayashi pseudo-metric  also  under such weaker regularity assumptions.\par
\smallskip
The structure of the paper is the following. 
 In \S 1, we recall  a few basic properties of  Monge-Amp\`ere exhaustions and Kobayashi metrics. In \S 2 we present our results 
with a short description of their proofs. Finally, in   \S 3 we present some open questions  that  might be 
 addressed using the results in \S 2. 
\par
\medskip
\section{Monge-Amp\`ere exhaustions and Kobayashi metrics}
\subsection{Monge-Amp\`ere exhaustions and Monge-Amp\`ere foliations} \label{sub2.1}
Let $D \subset \bC^n$ be a bounded domain with  boundary of class $\cC^{k,\a}$, $k \geq 2$, $\a> 0$.
A {\it Monge-Amp\`ere exhaustion}   for $\overline D$ is a continuous  exhaustion $\t: \overline D \to [0,1]$ satisfying the following conditions: 
\begin{itemize}
\item[i)]  The boundary $\p D$ coincides with the level set $ \{\t =1\} $  while the level set $\{\t = 0\}$ consists of exactly one interior point $z_o$, called {\it center}
of the exhaustion.
\item[ii)]  
The map $\t$ is of class $\cC^{k,\a}$ on $\overline D \setminus \{z_o\}$ and,  in general,  is only continuous at $z_o$. However, 
if $\pi: \wt{\overline D} \to \overline D$ denotes the blow up of  $\overline D$ at the center $z_o$, we  assume that the $\cC^{k,\a}$ lifted map 
$$\wt \t = \t|_{\overline D \setminus \{0\}} \circ \pi: \wt{\overline D}\setminus \{z_0\} \to  (0,1]$$
 admits a $\cC^{k,\a}$-extension to the whole $ \wt{\overline D}$. \\[-10pt]
\centerline{
\hskip-2.5cm 
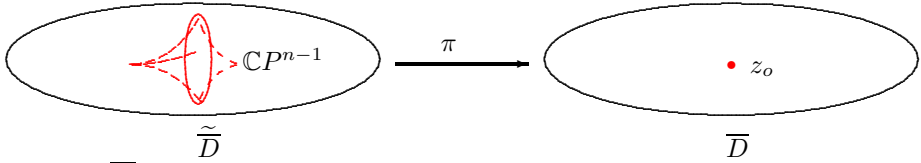
\begin{picture}(350,100)(0,0)
 \setlength{\unitlength}{1pt}
 \linethickness{0.5mm}
  \linethickness{0.3mm}
\thinlines
\put(75,40){ \renewcommand{\xscale}{2} 
\renewcommand{\xscaley}{-2} 
\renewcommand{\yscale}{0.6}
 \renewcommand{\yscalex}{0.6}
  \scaleput(10,10){\bigcircle{50}}
}
{\color{red}  \put(65,52){ \renewcommand{\xscale}{0.6} 
\renewcommand{\xscaley}{0.6} 
\renewcommand{\yscale}{2}
 \renewcommand{\yscalex}{-2}
  \scaleput(10,10){\bigcircle{12}}
}
  {\linethickness{0.1mm} \curvedashes[0.25mm]{2.5,5}
  \curve(50, 50, 65, 54, 77, 67)}
   {\linethickness{0.1mm} \curvedashes[0.25mm]{2,5}
  \curve(50, 50, 63, 51, 77, 55)}
   {\linethickness{0.1mm} \curvedashes[0.25mm]{2.5,5}
  \curve(50, 50, 65, 47, 78, 36)}
   {\linethickness{0.1mm} \curvedashes[0.2mm]{5,5}
  \curve(95, 50, 85, 54, 77, 67)}
   {\linethickness{0.1mm} \curvedashes[0.2mm]{5,5}
  \curve(95, 50, 85, 47, 77, 36)}
  }
%
  \put(275,40){\renewcommand{\xscale}{2} 
\renewcommand{\xscaley}{-2} 
\renewcommand{\yscale}{0.6}
 \renewcommand{\yscalex}{0.6}
  \scaleput(10,10){\bigcircle{50}}
  }
  \put(90, 48){\small $\bC P^{n-1}$}
   \put(73, 15){\small $\wt{\overline D}$}
  \put(273, 15){\small $\overline D$}
{\color{red}\put(275,50){\circle*{3}}}
\put(282, 47){\small $z_o$}
 \linethickness{0.4mm}
\put(165, 55){\small $\pi$}
\put(148, 50){\vector(1,0){50}}
 \end{picture}
 }
\ \\[- 1.5cm]
\item[iii)] On  $\overline D \setminus \{z_o\}$ the following differential conditions hold:
\begin{itemize}
\item[(a)] $2 i \p \bar \p \t > 0$; 
\item[(b)] $2 i \p \bar \p \log \t \geq 0$;  
\item[(c)] $(\p \bar \p \log \t)^n = 0$ ({\it homogeneous complex Monge-Amp\`ere equation}). 
\end{itemize}
Note that (a) - (c) imply that, at each point  $z$, the kernel of the $2$-form $ \p \bar \p \log \t_z$  is $1$-dimensional and that each  level set $\{\t = c\}$   is a strongly pseudoconvex real hypersurface.
\item[iv)] In proximity of the center,  $\log \t$ goes as $\log \t(z) \simeq \log(\| z - z_o\|) + O(1)$.
\end{itemize}
The closure  of a domains $D \subset \bC^n$, for which there is at least one Monge-Amp\`ere exhaustion,  is  called {\it (closed) domain of circular type}. \par
\begin{rem} The above definition, given here   for domains in $\bC^n$,  can be easily extended and  stated  in full generality  for complex manifolds with boundary. We refer to \cite{PS3} for details. 
\end{rem}
Up to a  few  minor  changes, the above notion of domain of circular type coincides with the one introduced by    the first author in 
\cite{Pt}  to capture the most crucial properties of the  following two important classes of domains. \par  
\medskip
\noindent{\bf Class A.} Let $\overline D \subset \bC^n$ be the closure of a complete circular domain,  i.e. of a domain $D = \{z:\mu(z) < 1\}$
 determined by  a defining  function $\mu: \bC^n \longrightarrow  [0, +\infty)$ satisfying the condition
 $$\mu(\lambda z) = |\lambda| \mu(z)\qquad \text{for all}\  \lambda \in \bC\ .$$
Assume  that the function  $\mu$,   called  the {\it Minkowski function} of $D$, has the following two properties:
\begin{itemize}
 \item[(a)]   one (and, consequently,   all) of the level sets $\{\mu = c\}$ for $0 < c \leq 1$  is a strongly pseudoconvex hypersurface,
 \item[(b)] $\mu$ is of class $\cC^{k,\a}$, with $k \geq 2$,  $\a > 0$,  on $\bC^n \setminus \{0\}$. 
 \end{itemize}
 One can then directly see that  the square  $\t \= (\mu|_{\overline D})^2$ is a Monge-Amp\`ere exhaustion for $\overline D$, centered  at  $z_o = 0$, so that      $\overline D$ is  of circular type.
\par
\medskip
\noindent{\bf Class B.}
Let $\overline D \subset \bC^n$ be  the closure of a strictly linearly convex domain $D$ with boundary of class $\cC^{k,\a}$ for some  $k \geq 4$,  $\a > 0$. Let also  $\d: D {\times} D \to \bR$ be the  Kobayashi pseudodistance of $D$ and, for each given $z_o \in D$, set 
\beq \label{MAz_o} \t^{(z_o)}: \overline D \longrightarrow [0, 1]\ ,\qquad \t^{(z_o)}(w) := \tanh^2(\d(z_o, w))\ . \eeq
Lempert's theory of  Kobayashi distance on convex domains (\cite{Le}) shows that {\it for each $z_o \in D$ the corresponding real  function \eqref{MAz_o} 
is a Monge-Amp\`ere exhaustion for  $\overline D$, centered at $z_o$}. Thus, any such domain is of circular type. Note  that, in contrast with the previous construction, for such domains  {\it any}  $z_o$ occurs  as the center of a  Monge-Amp\`ere exhaustion. 
\par
\medskip
An arbitrary  domain  $\overline D$ of circular type has always the following
crucial features. For any  $z \in \overline D \setminus \{z_o\}$,  let  $\cZ_z$  be the $1$-dimensional kernel
$$\cZ_z:= \ker(2 i \p \bar \p \log \t)|_z \subset T_z D\ .$$
Being formed by the  kernels of a closed $2$-form, the  distribution $\cZ$  is  integrable  and the (closures of) its integral leaves 
form a singular foliation of $D$, called {\it Monge-Amp\`ere foliation}. The leaves  are   
complex curves, each of them biholomorphic to the unit disk $\D$ in  $\bC$,    passing through the center $z_o$ of the Monge-Amp\`ere exhaustion  
(\cite{Pt, Pt1, PS}). \\
\centerline{
 \begin{picture}(350,80)(0,10)
 \setlength{\unitlength}{1pt}
 \linethickness{0.5mm}
  \linethickness{0.3mm}
\thinlines
  \put(80,40){\renewcommand{\xscale}{2} 
\renewcommand{\xscaley}{-2} 
\renewcommand{\yscale}{0.6}
 \renewcommand{\yscalex}{0.6}
  \scaleput(10,10){\bigcircle{50}}
  }
  {\color{red}\put(80,50){\circle*{5}}
\put(89, 50){\small $z_o = 0$}}
{  \linethickness{0.22mm}
 \curve(50,32, 80,50, 110,71)
 \curve(80,31, 80,50, 80,72)
 \curve(15,60, 80,50, 137,41)}
\put(205,15){
\curve(12,18, 30,10, 60, 12, 80, 14, 110, 25, 122, 36, 122, 42,  110, 53,
100, 58, 90, 61, 80, 62, 70, 63, 60, 64, 50, 63.5, 
47, 63,  44, 62, 40, 61, 35, 60,30, 57, 25, 53, 20, 49, 17, 46, 16, 45, 14, 42, 13, 41, 
10, 38, 9, 36, 8, 34, 8, 32, 8, 30, 8, 28,9,25, 10,23,11,20, 12,18)
}
 {\color{red}\put(240,40){\circle*{5}}
\put(250, 35){\small any $z_o$}}
{  \linethickness{0.22mm}
 \curve(230,26, 240,40, 250,77)
 \curve(252,26, 240,40, 216.5,53)
 \curve(217.5,32.5, 240,40, 328,55)
 }
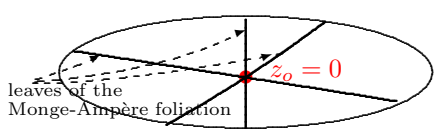
   \put(0, 15){\tiny  \bf strongly pseudoconvex circular domain}
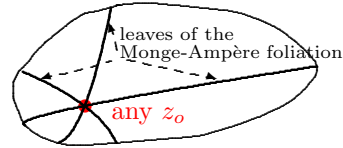
  \put(203, 15){\tiny  \bf strictly linearly convex domain}
  \put(-10, 43){\tiny leaves of the}
  \put(-10, 36){\tiny Monge-Amp\`ere foliation}
     {\linethickness{0.1mm} \curvedashes[0.2mm]{5,5}
  \curve(-3, 48, 15, 54, 25, 58)
   \curve(-3, 48, 35, 50, 94, 59)
     \curve(-3, 48, 45, 55, 80, 68)
  }
  \put(20,56){\vector(2,1){5}}
   \put(85,57){\vector(3,1){5}}
     \put(75,66){\vector(3,1){5}}
   \put(253, 65){\tiny leaves of the}
  \put(253, 58){\tiny Monge-Amp\`ere foliation}
  {\linethickness{0.1mm} \curvedashes[0.2mm]{5,5}
  \curve(253, 58, 240, 54, 222, 50)
  \curve(253, 58, 250, 70)
   \curve(253, 58, 290, 50)
    \put(228,52){\vector(-3,-1){5}}
 \put(250,69){\vector(-1,3){1}}
    \put(285,52){\vector(3,-1){5}}
  }
 \end{picture}
 }
The existence of such   foliation  has important consequences on the Kobayashi  pseudo-metrics of  domains of circular type, which we now  review.\par
\medskip 
\subsection{Kobayashi pseudo-metrics} 
Let $D$ be a bounded domain of $\bC^n$ and $z_o$ a point of $D$. 
We recall that the {\it Kobayashi (infinitesimal) pseudo-metric} of $D$ at $z_o$  is the real valued function on $T_{z_o} D \simeq \bC^n$, defined  by
$$\k_{z_o} : T_{z_o} D \setminus \{\underline 0\}  \longrightarrow [0, + \infty) \ ,\qquad  \k_{z_o}(\underline v) := \inf_{f \in \cA_{(z_o, \underline v)}} \|f^{-1}_*(\underline v)\|\ ,$$
where  $\cA_{(z_o, \underline v)}$ denotes the set  of all holomorphic maps from $\D$ to $D$, passing through $z_o$ and tangent to $\underline v$, i.e. 
\beq\nonumber
\cA_{(z_o, \underline v)} \!\!:{=}\bigg\{f: \D \to D \ \text{ holom. with }  f(0) = z_o\ ,\ f_*\left(\l \frac{\p}{\p x}\bigg|_0\right) = \underline v\ \text{for}\ \l \in \bR\ \bigg\}\ .
\eeq
The literature on the Kobayashi pseudometric is vast (see e.g.   \cite{Ko, Ko1, JP} and references therein). Here we just recall  two   simple -- but  crucial -- facts. \par
\smallskip
\begin{itemize}[itemsep=8pt plus 5pt minus 2pt, leftmargin=18pt]
\item[1)] It satisfies the so-called {\it distance decreasing property}, i.e. for any holomorphic map  $F:D \longrightarrow D'$  between  two   domains $D, D'$ and  for each $z_o \in D$, the Kobayashi pseudo-metrics $\k_{z_o}$ and $\k_{F(z_o)}$ of $D$ and $D'$, respectively,  satisfy the inequality 
 $$\k_{F(z_o)}(F_*(\underline v)) \leq \k_{z_o}(\underline v)\qquad \text{for each}\ \ \underline v \in T_{z_o} D\ .$$
 It follows immediately that if  $F$  is   a biholomorphism, then the  equality holds, meaning that the Kobayashi pseudo-metric is a (very important) biholomorphic invariant.
\item[2)] For  each $z_o \in D$ and $\underline v \in T_{z_o} D$, 
$$\k_{z_o}(\l \underline v) =\! |\l| \,\k_{z_o}(\underline v)\qquad \text{for each}\ \l \in \bC\ .$$
This  yields that  the {\it indicatrix} of the Kobayashi pseudo-metric at   $z_o$, that is the set 
 \beq\label{indicatrices} I_{z_o} := \{ \underline v \in T_{z_o} D \setminus \{0\}\ : \k_{z_o}(\underline v) < 1\} \cup\{\underline 0\}\ ,\eeq
is always a balanced domain of $T_{z_o} D \simeq \bC^n$.  
\end{itemize}
The situations where  $I_{z_o}$ is a {\it strongly pseudoconvex  domain} for each point $z_o$ are of particular interest, because in those cases 
 the Kobayashi pseudo-metric $\k: TD \to [0, + \infty)$ is a (possibly non-smooth) complex Finsler metric.\par
\smallskip
\subsection{Finsler metrics}
Let us shortly recall   the definitions of real and complex Finsler metrics.  Given an $n$-dimensional real manifold  $M$, with tangent bundle $TM$,   a  (smooth) {\it  real Finsler metric}   is a continuous map 
$$F: TM \longrightarrow [0, + \infty)\ ,$$
which is   $\cC^\infty$  on $TM^o\= TM \setminus \{\text{zero section}\}$  and  such that
\begin{itemize}
\item[a)] for every non-negative {\it real} number  $\l$ and every   vector $(x, \underline v) \in T_{x} M$  at some point $x \in M$ one has that 
$$F(x; \l \underline v) = |\l| F(x;  \underline v)$$ 
\item[b)] for each  $x_o \in M$, the indicatrix $ I_{x_o}\! :{=}\{ (x_o, \underline v): F(x_o; \underline v) < 1\} $ is a strictly linearly convex domain of $T_{x_o} M \simeq \bR^n$. 
\end{itemize}
In other words,  a smooth real Finsler metric  is a norm function on all tangent spaces of $M$, which is   smoothly depending on the base points and on the (non-zero) tangent vectors (when such  dependence is not  $\cC^\infty$, it is usually  said  that  $F$ is  a  ``non-smooth'' Finsler metric). In fact, a very simple example of Finsler metric on a manifold $M$ is  given by the norm function 
$$F(x; \underline v) \= \sqrt{g_{x}(\underline v, \underline v)}\ ,$$
determined by some fixed Riemannian  metric $g$ on $M$. But many other examples, for which there is no associated Riemannian metric,  can be easily  constructed. Indeed, in order to define a Finsler metric, it suffices to fix a   linearly convex indicatrix $I_x$ in each tangent space, smoothly depending on the base point $x$, and use it to define a norm function. If the assigned indicatrices are not linearly equivalent to quadrics, the corresponding Finsler metric cannot be associated with any Riemannian metric.  \par
\smallskip
The notion of complex Finsler metric is very similar. If $M$ is a complex manifold of complex dimension $n$, a  (smooth) {\it  complex Finsler metric} on $M$ is a continuous  real valued map 
$$F: TM \longrightarrow [0, + \infty)\ ,$$
which is  $\cC^\infty$  on $TM^o\= TM \setminus \{\text{zero section}\}$  and satisfies the following   analogues of (a) and (b): 
\begin{itemize}
\item[a')] for every   {\it complex} number  $\l$ and every  vector $(x, \underline v) \in T_x M$, $x \in M$, 
$$F(x; \l \underline v) = |\l| F(x; \underline v);$$ 
\item[b')] for each  $x_o \in M$, the indicatrix $ I_{x_o} := \{ (x_o, \underline v): F(x_o; \underline v) < 1\}$ is a strongly pseudoconvex  domain  of $T_{x_o} M \simeq \bC^n$. 
\end{itemize}
As before, a very simple example of complex Finsler metric on a complex manifold $M$ is  given by the norm function 
$$F(x; \underline v) \= \sqrt{h_{x}(\underline v, \underline v)}\ ,$$
determined by an Hermitian  metric $h$ on $M$. But many other examples can be easily  determined, for which there is no associated Hermitian metric.
In perfect analogy with the real case, complex Finsler metrics are uniquely determined by  the associated family of indicatrices.  \par
\smallskip
As for the  Kobayashi metric, also the literature on Finsler metrics is  enormous. For an introduction to this important and interesting  area,  the  reader might take a look at   standard texts as  \cite{Ch, BCS, AP}  and references therein. \par
\medskip
\subsection{Kobayashi metrics and Monge-Amp\`ere foliations}
Coming back to the Kobayashi pseudo-metric $\k$  of a domain $D \subset \bC^n$, if the  
indicatrices \eqref{indicatrices} are  smooth and strongly pseudoconvex and if  they smoothly depend on their base points, then $\k$ is a complex Finsler metric. But here come two of the most  unfriendly features  of Kobayashi pseudo-metrics.
\begin{itemize}
\item[--] For a generic domain, the indicatrix  \eqref{indicatrices}  is usually {\it not strongly pseudoconvex} and  the function  $\k: T D^o = TD \setminus \{\text{zero section}\} \to [0, + \infty)$  is  often no better than $\cC^0$.
\item[--] Domains, for  which  $\k$ can be explicitly  computed -- hence   analyzed  in greater  detail -- are not easy to find.
\end{itemize}
Nonetheless, for  domains of circular type things are somehow nicer. 
\par
\medskip
Let $\overline D \subset \bC^n$ be a strongly pseudoconvex domain of circular type, equipped with a Monge-Amp\`ere exhaustion  $\t$ with center  $z_o$. As it is shown in  \cite{Pt1},  
 for each  non-zero tangent vector $\underline v \in T_{z_o} D$  at the center $z_o$, there exists a {\it unique} proper holomorphic disk  
$f^{(\underline v)}{:} \overline \D {\to} \overline D$,
 whose image $f^{(\underline v)}(\overline \D)$  coincides with  the closure of a leaf of the Monge-Amp\`ere foliation of  $\t$ and tangent to $\underline v$ at the origin, i.e.  satisfying the conditions
 \beq \label{stat} f(0) = z_o\ , \qquad f_*\left(\l\frac{\p}{\p x}\bigg|_{0}\right) = \underline v\qquad \text{for some}\ \l \in \bR\ .\eeq \\
Here  $\l$ is  uniquely determined by $\underline v$ and  is  non-zero. Let us  denote it by $\l^{(\underline v)}$. 
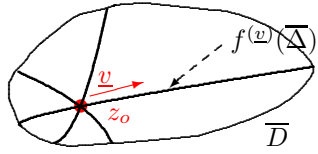
\par
\centerline{
 \begin{picture}(140,75)(0,10)
\put(5,15){
\curve(12,18, 30,10, 60, 12, 80, 14, 110, 25, 122, 36, 122, 42,  110, 53,
100, 58, 90, 61, 80, 62, 70, 63, 60, 64, 50, 63.5, 
47, 63,  44, 62, 40, 61, 35, 60,30, 57, 25, 53, 20, 49, 17, 46, 16, 45, 14, 42, 13, 41, 
10, 38, 9, 36, 8, 34, 8, 32, 8, 30, 8, 28,9,25, 10,23,11,20, 12,18)
}
 {\color{red}\put(40,40){\circle*{5}}
 \put(50, 35){\small  $z_o$}}
 \put(110, 25){\small $\overline D$}
{  \linethickness{0.22mm}
 \curve(30,26, 40,40, 50,77)
 \curve(52,26, 40,40, 16.5,53)
 \curve(17.5,32.5, 40,40, 128,55)
 }
 {\color{red} \put(40,44){\vector(4,1){20}}
  \put(43, 49){\small  $\underline v$}}
  \put(93, 62){\small $f^{(\underline v)}(\overline \D)$}
   {\linethickness{0.1mm} \curvedashes[0.2mm]{5,5}
  \curve(92, 65, 73, 50)
 \put(75,52){\vector(-1,-1){5}}
  }
 \end{picture}
 }
\noindent  A crucial relation between the Monge-Amp\`ere foliation  of $\overline D$ and its Kobayashi pseudo-metric is represented by the following  
two facts (\cite{Pt3}): 
\begin{itemize}
\item[a)]  for each $ \underline v \in T_{z_o} D \setminus \{0\}$ one has that 
$\k_{z_o}(\underline v) = \l^{(\underline v)}$; 
\item[b)]  the indicatrix $I_{z_o}$  at the center of $\k$  is  strongly pseudoconvex.
\end{itemize}
\par
\smallskip
\noindent This   yields to the following couple of nice properties.
\begin{itemize}[itemsep=4pt plus 5pt minus 2pt, leftmargin=18pt]
\item[--] {\it If $D \subset \bC^n$ is a smoothly bounded, strongly pseudoconvex    circular domain, then  at $z_o = 0$ the indicatrix $I_{z_o = 0}$  is smooth and  strongly pseudoconvex. In fact, it is linearly equivalent to the domain $D$. }
\item[--] {\it  If $D \subset \bC^n$ is  a  smoothly bounded strictly  linearly  convex   domain,  then at each point $z_o \in D$,    the indicatrix $I_{z_o}$  is smooth and strongly pseudoconvex and it smoothly depends  on the base point. In  other words, $\k$  is a smooth  complex Finsler metric. }
\end{itemize}
\par
\medskip
\section{The phenomenon of propagation of regularity}
\subsection{Domains with a lot of Monge-Amp\`ere exhaustions} In this section, we study the properties of the 
domains of circular type admitting a one-parameter family of  Monge-Amp\`ere exhaustions,   centered at the   points of a given curve. We  will shortly see 
that the existence of such Monge-Amp\`ere exhaustions is equivalent to the existence of a special one-parameter family of homeomorphisms of the domain into 
itself, with  nice regularity properties outside the curve  and satisfying a special set   of differential constraints. On the basis of these facts,  all main  results presented in this note will be built. 
\subsubsection{Domains with two Monge-Amp\`ere exhaustions}
 Let   $\overline D \subset \bC^n$ be a domain of circular type and assume that  it admits   at least  {\it two} distinct  Monge-Amp\`ere exhaustions,  the first centered at $z_o$ and the second at $w_o$, and consequently  with  {\it two}  Monge-Amp\`ere  foliations, one for $z_o$, the other for $w_o$.\\
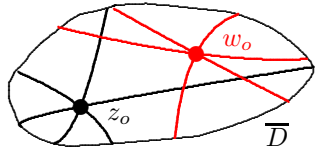
\centerline{
 \begin{picture}(140,75)(0,10)
\put(5,15){
\curve(12,18, 30,10, 60, 12, 80, 14, 110, 25, 122, 36, 122, 42,  110, 53,
100, 58, 90, 61, 80, 62, 70, 63, 60, 64, 50, 63.5, 
47, 63,  44, 62, 40, 61, 35, 60,30, 57, 25, 53, 20, 49, 17, 46, 16, 45, 14, 42, 13, 41, 
10, 38, 9, 36, 8, 34, 8, 32, 8, 30, 8, 28,9,25, 10,23,11,20, 12,18)
}
\put(40,40){\circle*{6}}
 \put(50, 35){\small  $z_o$}
 \put(110, 25){\small $\overline D$}
{  \linethickness{0.22mm}
 \curve(30,26, 40,40, 50,77)
 \curve(52,26, 40,40, 16.5,53)
 \curve(17.5,32.5, 40,40, 128,55)
 }
{\color{red} \put(80,60){\circle*{6}}
 \put(90, 63){\small  $w_o$}
{  \linethickness{0.22mm}
 \curve(115,42, 80,60, 49,77)
 \curve(72,28, 80,60, 96,75)
 \curve(28,70, 80,60, 122,58)
 }
 }
 \end{picture}
 }
Applying various  results by the first author on the so-called {\it circular representation} (\cite{Pt1}), it is possible to show  that  the existence of  a  Monge-Amp\`ere exhaustion $\t^{(z_o)}$ of class $\cC^{k, \a}$,  $k \geq 4$, $\a> 0$,   off the center $z_o$, implies also   the existence  of a    ``straightening'' homeomorphism $\Psi^{(z_o)}: \overline D \longrightarrow \overline \bB^n$  onto the closed unit ball of  $\bC^n$  centered at  $0$, 
with the following nice properties (\cite{PS, PS3}): 
\begin{itemize}[itemsep=4pt plus 5pt minus 2pt, leftmargin=18pt]
\item[a)] $\Psi^{(z_o)}(z_o) = 0$ and each   level set $\{\t^{(z_o)}(z) = c\}$, $0 < c \leq 1$,  is mapped diffeomorphically  onto the  sphere $\{|z| = c \}$ of radius $c$. The restriction of $\Psi^{(z)}$ on each such  level set is  {\it in general not a CR map}, but nonetheless maps the real   contact distributions underlying the two CR structures one into the other.
\item[b)]  The disks of the Monge-Amp\`ere foliation through $z_o$  are mapped bijectively onto the straight disks of $\bB^n$  through  $0$. Each restriction of $\Psi^{(z_o)}$ along one such a  disk  is a biholomorphism.
\item[c)]   $\Psi^{(z_o)}$ is of class $\cC^{k-2, \a}$ on $\overline D \setminus \{z_o\}$. Further, if we denote by  $\wt{\overline D}^{(z_o)}$ and $\wt{\overline \bB^n}$  the blow-ups at $z_o$ and $0$ of the  domains, then the restriction $\Psi^{(z_o)}|_{\overline D \setminus \{z_o\}}$ admits a unique  $\cC^{k-2, \a}$-extension to a map  between $\wt{\overline D}^{(z_o)}$ and $\wt{\overline \bB^n}$.
\end{itemize}
\ \\[-1.9cm]
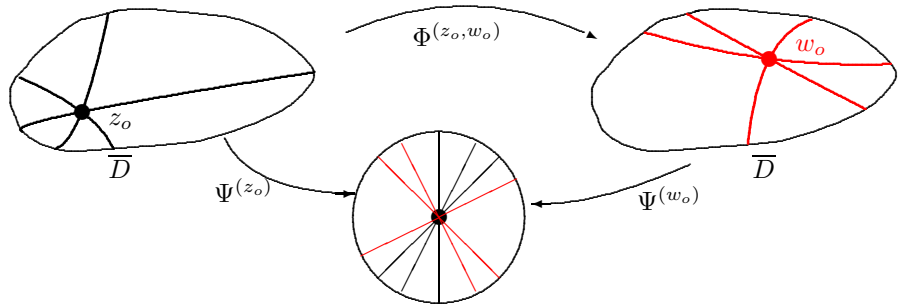
\centerline{
 \begin{picture}(360,175)(0,60)
\put(5,115){
\curve(12,18, 30,10, 60, 12, 80, 14, 110, 25, 122, 36, 122, 42,  110, 53,
100, 58, 90, 61, 80, 62, 70, 63, 60, 64, 50, 63.5, 
47, 63,  44, 62, 40, 61, 35, 60,30, 57, 25, 53, 20, 49, 17, 46, 16, 45, 14, 42, 13, 41, 
10, 38, 9, 36, 8, 34, 8, 32, 8, 30, 8, 28,9,25, 10,23,11,20, 12,18)
\put(35,25){\circle*{6}}
 \put(45, 20){\small  $z_o$}
 \put(45, 0){\small $\overline D$}
{  \linethickness{0.22mm}
 \curve(25,11, 35,25, 45,62)
 \curve(47,11, 35,25, 11.5,38)
 \curve(12.5,17.5, 35,25, 123,40)
 }
}
\put(225,115){
\curve(12,18, 30,10, 60, 12, 80, 14, 110, 25, 122, 36, 122, 42,  110, 53,
100, 58, 90, 61, 80, 62, 70, 63, 60, 64, 50, 63.5, 
47, 63,  44, 62, 40, 61, 35, 60,30, 57, 25, 53, 20, 49, 17, 46, 16, 45, 14, 42, 13, 41, 
10, 38, 9, 36, 8, 34, 8, 32, 8, 30, 8, 28,9,25, 10,23,11,20, 12,18)
{\color{red} \put(75,45){\circle*{6}}
 \put(85, 48){\small  $w_o$}
{  \linethickness{0.22mm}
\curve(111.5,26, 75,45, 44,62)
\curve(67,13, 75,45, 91,60)
 \curve(28,55, 75,45, 121,43)
 }
 }
  \put(65, 0){\small $\overline D$}
  }
  \put(175,100){\bigcircle{65}}
   \put(175,100){\circle*{6}}
\put(175,100){\line(-1,-1){23}}
\put(175,100){\line(1,1){23}}
\put(175,100){\line(0,-1){32.5}}
\put(175,100){\line(0,1){32,5}}
\put(175,100){\line(-1,-2){14}}
\put(175,100){\line(1,2){14}}
{\color{red}
\put(175,100){\line(-2,-1){29.5}}
\put(175,100){\line(2,1){29.5}}
\put(175,100){\line(1,-1){23}}
\put(175,100){\line(-1,1){23}}
\put(175,100){\line(1,-2){14}}
\put(175,100){\line(-1,2){14}}
}
 \curve(140,170, 185,180, 230, 170)
 \put(228, 171){\vector(2,-1){7}}
  \curve(270,120, 240,109,215, 105)
 \put(217, 105){\vector(-1,0){7}}
  \curve(94,130, 110,115,140, 109)
  \put(136, 109){\vector(1,0){7}}
  \put(165, 165){\small $\Phi^{(z_o, w_o)}$}
  \put(90, 106){\small $\Psi^{(z_o)}$}
  \put(250, 103){\small $\Psi^{(w_o)}$}
 \end{picture}
 }
Due to this, if we have   {\it two} Monge-Amp\`ere exhaustions, we may   compose  the {\it two} associated ``straightening'' homeomorphisms and get
a homeomorphism from $\overline D$ into itself
$\Phi^{(z_o, w_o) }:= (\Psi^{(w_o)})^{-1} \circ \Psi^{(z_o)}: \overline D \longrightarrow \overline D$
such that 
\begin{itemize}[itemsep=4pt plus 5pt minus 2pt, leftmargin=18pt]
\item[$\a$)] $\Phi^{(z_o, w_o)}$ transforms the center $z_o$ into the center $w_o$ and maps  each   level set of $\t^{(z_o)}$ into the corresponding one of 
$\t^{(w_o)}$. The restriction of $\Phi^{(z_o, w_o)}$ on each level set maps one to the other the contact structures underlying the CR structures of those hypersurfaces. 
\item[$\b$)]  The disks of the Monge-Amp\`ere foliation  through $z_o$   are mapped bijectively onto the disks of the Monge-Amp\`ere foliation  through $w_o$
and on each of them the restriction of $\Phi^{(z_o, w_o) }$  is a biholomorphism. 
\item[$\g$)] the map $\Phi^{(z_o, w_o)}$ is of class $\cC^{k-2, \a}$ on $\overline D \setminus \{z_o\}$ and, considering the blow-ups $\wt{\overline D}^{(z_o)}$ and $\wt{\overline D}^{(w_o)}$ at $z_o$ and $w_o$, the restriction $\Phi^{(z_o, w_o)}|_{\overline D \setminus \{z_o\}}$ admits a unique $\cC^{k-2, \a}$-extension that  goes from  $\wt{\overline D}^{(z_o)}$  to  $\wt{\overline D}^{(w_o)}$.
\end{itemize}

\par
\medskip
\subsubsection{Domains with a one-parameter family of  Monge-Amp\`ere exhaustions}\label{sect312}
 Let us now address  the case of a closed domain of circular type $\overline D$ admitting a  {\it one-parameter family of Monge-Amp\`ere exhaustions}, centered at the points $z_t$, $t \in [0,1]$, of a smooth curve of $D$.  \\
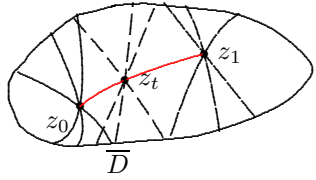
\centerline{
 \begin{picture}(125,80)(0,0)
\put(-2,3){
\curve(12,18, 30,10, 60, 12, 80, 14, 110, 25, 122, 36, 122, 42,  110, 53,
100, 58, 90, 61, 80, 62, 70, 63, 60, 64, 50, 63.5, 
47, 63,  44, 62, 40, 61, 35, 60,30, 57, 25, 53, 20, 49, 17, 46, 16, 45, 14, 42, 13, 41, 
10, 38, 9, 36, 8, 34, 8, 32, 8, 30, 8, 28,9,25, 10,23,11,20, 12,18)
\put(35,25){\circle*{3}}
 \put(22, 17){\small  $z_0$}
 \put(45, 0){\small $\overline D$}
{  \linethickness{0.21mm}
 \curve(25,11, 35,25, 34,60)
 \curve(47,11, 35,25, 11.5,38)
\curve(34,10, 35,25, 24,53)
 }
 \put(52,35){\circle*{3}}
 \put(57, 33){\small  $z_t$}
 \put(82,45){\circle*{3}}
 \put(87, 43){\small  $z_1$}
{\color{red}
 \curve(35,25,52,35,82,45)
 }
 {\linethickness{0.2mm} \curvedashes[0.5mm]{2,5}
  \curve(40, 8, 52,35, 63, 62)
\curve(48, 9, 52,35, 55, 64)
\curve(68, 11, 52,35, 33, 58)
}
 { \linethickness{0.2mm} \curvedashes[0.5mm]{1,5}
  \curve(69, 13, 82,45, 95, 60)
\curve(88, 14.5, 82,45, 75, 63)
\curve(102, 20, 82,45, 67, 64)
}
}
 \end{picture}
 }
The   remarks of previous section imply that  in this case  $\overline D$  is equipped with  a $1$-parameter family of homeo\-morphisms 
$\Phi^{(t)}: \overline D \longrightarrow \overline D$,  $t \in [0,1]$,  such that:
\begin{itemize}[itemsep=4pt plus 5pt minus 2pt, leftmargin=18pt]
\item[1)]   $\Phi^{(t)}$ is  $\cC^{k-2, \a}$ on  $\overline D \setminus \{z_0\}$, with    $\Phi^{(t)}|_{\overline D \setminus \{z_0\}}$ with a  $\cC^{k-2, \a}$-extension   to a map from the blow up  $\wt{\overline D}^{(z_0)}$  at $z_0$ onto the blow up  $\wt{\overline D}^{(z_t)}$  at $z_t$; 
\item[2)] $\Phi^{(t)}$ maps $z_0$ into $z_t$, and  sends the   level sets of $\t^{(z_o)}$ onto the corresponding level sets of 
$\t^{(w_o)}$ by contact transformations. In particular, it induces a $\cC^{k-2, \a}$ contact map from  $\p D$ into itself and  $\t^{(t)} = \t^{(0)}\circ (\Phi^{(t)})^{-1}$. 
\item[3)]  The disks of the  Monge-Amp\`ere foliation   through $z_0$ are biholomorphically mapped  by $\Phi^{(t)}$    onto the disks of the Monge-Amp\`ere foliation  through $z_t$. 
\end{itemize}
 Since each  of the associated lifts   between  blow-ups  is $\cC^{k-2, \a}$ with $k-2 \geq 2$,     we may consider  the pull-backed tensor fields on $\wt{ \overline D}^{(z_0)}$ defined by 
\beq \label{Jt} J_t  \= (\Phi^{(t)-1})_*(J_o)\ , \qquad t \in [0,1]\ ,\eeq
where  $J_o$ stands for the standard complex structure of the blow up $\wt{\overline D}^{(z_t)}$  at $z_t$.  
The $J_t$   are tensor fields of type $(1,1)$,   they verify the condition  $J^2_t|_z = -I$ and their 
Nijenhueis tensors  are identically vanishing, being  each $J_t$ a pull-back of  the {\it integrable} complex structure  $J_o$. Further, each of them is  of class $\cC^{k-1, \a}$ with $k-1 \geq 1$ and $\a > 0$.  Hence, by   Newlander-Nirenberg Theorem  (\cite{NN, Ma,NW,We,HT}), each $J_t$  is  a {\it non-standard} integrable complex structure.  We call them   {\it non-standard} simply because in general  the  $\Phi^{(t)}$   {\it are not  biholomorphisms}  and, consequently,   the pull-backs of the standard complex $J_o$  by such maps are different from $J_o$. \par 
\medskip
If the curve of centers $z_t$  is  at least $\cC^1$ and if we select the diffeomorphisms $\Phi^{(t)}$ in such a way that the  family $\Phi^{(t)}$ is differentiable with respect to the parameter $t$, we may also consider the one-parameter  family of  vector field on $\overline D {\setminus} \{z_0\}$, defined by 
\beq \label{Xt} X_t|_z := (\Phi^{(t)-1})_*\bigg( \frac{d \Phi^{(t)}}{dt} \bigg|_{(t, z)}\bigg) \qquad  \text{for each}\ \ z \in \overline D {\setminus} \{z_0\}\ .\eeq
For getting a physical  intuition of such vector fields,   consider  the one-parameter  family of maps $\Phi^{(t)}$ as   a fluid motion and  the coordinates of  the points  of   $\overline D$   as  Lagrange coordinates for the fluid. In this way 
 $X_t$ can be interpreted as  the  {\it velocity field of the flow at time $t$ in Lagrangian coordinates}. 
\par
\smallskip
We will shortly see that  all crucial information about the  $\Phi^{(t)}$ is encoded in the  one-parameter family of pairs $(X_t, J_t)$, $t \in [0,1]$, which we   call the  {\it fundamental pair}  for the  one-parameter family  $\Phi^{(t)}$.
\par 
\medskip

\subsubsection{A simple  model example.} Assume that $\overline D = \overline \bB^n$ ,  let $v_o \neq 0$ in $\bB^n$ and denote by   $z_t = t{\cdot} v_o$, $t \in [0,1]$,  the points of the segment between $0$ and $v_o$. Since $\bB^n$ is homogeneous, we may  consider a smooth family of automorphisms $F^{(t)} \in \Aut(\overline{\bB^n}, J_o)$,   mapping  the origin into the   points $z_t$.\\
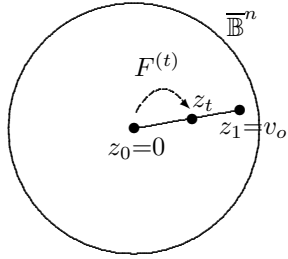
\centerline{
\begin{picture}
(100,100)(0,0)
\put(50,50){\bigcircle{95}}
 \put(50,50){\circle*{4}}
 \put(40,40){\small $z_0{=}0$}
  \put(72,53.5){\circle*{4}}
 \put(72,58){\small $z_t$}
 \put(90,57){\circle*{4}}
\put(82, 48){\small $z_1{=} v_o$}
 \curve(50,50, 90, 57)
 {\linethickness{0.2mm} \curvedashes[0.2mm]{1,5}
\curve(50,55, 60, 65, 70, 58.5)}
\put(71, 57.5){\vector(1,-1){0.2}}
\put(50, 70){\small $F^{(t)}$}
\put(85, 85){\small $\overline \bB^n$}
\end{picture}
}
$\overline \bB^n$ is clearly a circular domain with Minkowski function $\|\cdot\|$ and
$$\t^{(0)}: \overline \bB^n \longrightarrow [0,1]\ ,\qquad \t^{(0)}(w):=\|w\|^2$$
 is a Monge-Amp\`ere exhaustion for $\bB^n$  centered at  $0$. The corresponding   Monge-Amp\`ere  foliation is made by  the 
straight radial disks $f^{(\underline v)}(\overline\D)$,  images of the  maps $f^{(\underline v)}(\z) {:=}\z{\cdot} \underline v$ with $\|\underline v\| = 1$.
Using  the fact that each $F^{(t)}$ is a biholomorphism from $\bB^n$ into itself, one can directly check that  the exhaustions
$$\t^{(t)} = \t^{(0)} \circ F^{(t)-1}$$
are all  Monge-Amp\`ere exhaustions, each of them  centered at  a different  $z_t$.  The  corresponding  Monge-Amp\`ere  foliations are  made of the images under the maps $F^{(t)}$ of  the straight radial disks  through the origin.\\
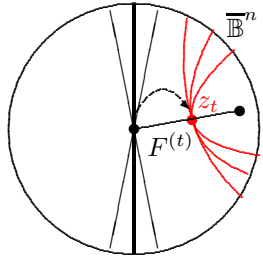
\centerline{
\begin{picture}
(100,100)(0,0)
\put(50,50){\bigcircle{95}}
 \put(50,50){\circle*{4}}
  {\color{red}
  \put(72,53.5){\circle*{4}}
 \put(74,58){\small $z_t$}
 \curve(70, 92, 72, 53.5, 91, 24)
 \curve(80, 87.5, 72, 53.5, 93, 33)
  \curve(88, 80, 72, 53.5, 97.5, 40)
 }
{  \linethickness{0.3mm}
\put(50,50){\line(0,1){47}}
\put(50,50){\line(0,-1){47}}
\put(50,50){\line(1,5){9}}
\put(50,50){\line(-1,-5){9}}
\put(50,50){\line(-1,5){9}}
\put(50,50){\line(1,-5){9}}
 }
  \put(90,57){\circle*{4}}
 \curve(50,50, 90, 57)
 {\linethickness{0.2mm} \curvedashes[0.2mm]{1,5}
\curve(50,55, 60, 65, 70, 58.5)}
\put(71, 57.5){\vector(1,-1){0.2}}
\put(55, 40){\small $F^{(t)}$}
\put(85, 85){\small $\overline \bB^n$}
\end{picture}
}
We are   in the situation considered in  the previous section: There is  a one-parameter family of Monge-Amp\`ere exhaustions $\t^{(t)}$, centered at the points of a curve $z_t$. In this  special case, we may take as  maps $\Phi^{(t)}: \overline \bB^n \to \overline \bB^n$    the biholomorphisms $\Phi^{(t)} = F^{(t)}$ and obtain as  fundamental pair
$$ J_t = (F^{(t)-1})_*(J_o) = J_o\ ,\qquad X_t = (F^{(t)-1})_*\left(\frac{d F^{(t)}}{dt}\right)\ .$$
This situation is however somehow  peculiar, because, in contrast with the generic case,  here each map $\Phi^{(t)} = F^{(t)}$ is a biholomorphism. 
\par
\medskip
\subsection{The Propagation of Regularity Theorem} We are now ready to state the main result and provide a  short description of the tools which are used to obtain it. Basically,  all  is built upon  the next two lemmas. 
\subsubsection{The  main  lemmas} If  $\Phi^{(t)}{:}\overline D {\to} \overline D$ is one of the  families of homeomorphisms considered in \S \ref{sect312},  the associated  fundamental pair $(X_t, J_t)$, $t \in [0,1]$,  satisfies  an important set of  differential constraints, which we are now going to  list.  In order to prevent    diversions of  reader's attention from   the most crucial aspects,  we give here just    intuitive descriptions of such constraints. The interested reader is referred to \cite{PS3} for  the detailed expressions. \par 
The constraints  that  any fundamental pair  satisfies are: 
\begin{itemize}[itemsep=4pt plus 5pt minus 2pt, leftmargin=18pt]
\item[(A)]  The restrictions of each   $J_t$ to  the tangent spaces of the leaves of the   Monge-Amp\`ere foliations of $\t^{(0)}$  coincide with the restrictions to the same spaces of the standard complex structure $J_o$. This is due to the fact that  the restriction of  $\Phi^{(t)}$  on  each  such a leaf is a biholomorphism.
\item[(B)] The one-parameter families $X_t$ and $J_t$ satisfy  $\frac{d J_t}{dt} =  \cL_{X_t} J_t$ with $J_{t = 0} = J_o$.
\item[(C)] Each $J_t$ has  identically vanishing Nijenhueis tensor, as  pointed  in \S \ref{sect312}. 
\item[(D)] Each vector field $X_t$ satisfies  special  boundary conditions at  $\p D$ and on its limit behavior  near $z_0$. They  are  determined 
by  the fact that each   $\Phi^{(t)}$ maps diffeomorphically $\p D$ into itself   and, at the same time, it admits  a  $\cC^{k-2,\a}$ extension between  the blow ups at $z_0$ and $z_t$.
\end{itemize}
 The  explicit expression for the  limit behavior of $X_t$  at $z_0$  mentioned in (D)  depends in a non trivial way {\it on the tangent vector of the curve $z_t$ at the time $t$}.  Such  dependence is quite technical and we  refer to \cite{PS3} for explicit details. \par
 \smallskip
 All this motivates the next
\begin{definition} 
Let $z_t$,  $t \in [0,1]$,  be   a $\cC^1$ curve $z_t$  in $D$ starting from the center $z_0$.  We call  {\it abstract fundamental pair  guided by $z_t$} any
 one-parameter family of pairs $(X_t, J_t)$,   formed by vector fields on $\overline D \setminus \{z_0\}$  and 
 almost complex structures $J_t$ on $\wt{\overline D}^{(z_0)}$    satisfying the  constraints (A) -- (D). \par
 If there exists  a one-parameter family $\t^{(t)}$ of Monge-Amp\`ere exhaustions, centered at the points $z_t$ and  with associated maps  $\Phi^{(t)}$ having  $(X_t, J_t)$ as fundamental  pair, we call the pair  a  {\it concrete fundamental pair}.
\end{definition} 
We may now state our two main  lemmas (\cite{PS3}).\par
\begin{lem} \label{mainLemma1} Let $\t = \t^{(0)}$ be a Monge-Amp\`ere  exhaustion on $\overline D$,  centered at  $z_0$ and of class   $\cC^{k, \a}$, $k \geq 4$, $\a > 0$,  on   $\overline D \setminus\{z_0\}$. 
Any  {\rm abstract} fundamental pair $(X_t, J_t)$ of class $\cC^{k-2,\a}$ is    {\rm concrete}  and is associated  with  the maps $\Phi^{(t)}$ of a one-parameter family of  Monge-Amp\`ere exhaustions $\t^{(t)}$ of class $\cC^{k-2,\a}$ off the centers.
\end{lem}
\begin{lem}  \label{mainLemma2}  Let $\t = \t^{(0)}$ be a Monge-Amp\`ere  exhaustion on $\overline D$ as in the previous lemma and denote by   $z_t$, $t \in [0,1]$, the points of a radious of a   holomorphic disk $f(\D)$  of the  Monge-Amp\`ere foliation   through   $z_0$  (i.e. $z_t$ has the form $z_t = f(tv)$  for some fixed $v$ with $\|v \| = 1$).
 Then there exists a value $\l_o \in (0, 1]$ such that, for all $\l \in (0, \l_o)$ there is an {\rm abstract} (hence  {\rm concrete}) fundamental pair $(X_t, J_t)$ of class   $\cC^{k-2, \a}$   guided by the curve $z_{\l t}$.
\end{lem}
The proof of Lemma \ref{mainLemma1} essentially consists   of two parts.  One first shows that,  for  a given  abstract pair $(X_t, J_t)$,  one can solve  the differential problem in $\Phi^{(t)}$, given by  the equation \eqref{Xt}  and the initial conditions $\Phi^{(t = 0)} = \Id$. For proving this, the key idea is to observe that the  {\it non-linear} equation \eqref{Xt} is  equivalent to a  {\it quasi-linear}  equation on  the inverse maps  $\Psi^{(t)} = (\Phi^{(t)})^{-1}$, for which the existence of solutions can be  proved with little effort.  Secondly,  one tries to show that the compositions  $\t^{(t)} \= \t^{(0)} \circ (\Phi^{(t)})^{-1}$ satisfy all conditions for being 
Monge-Amp\`ere exhaustions.  The only non immediate points of this check are reduced  to prove that the maps $\Phi^{(t)}$ have uniformly bounded  Jacobians at the points where they are differentiable.  This is   first proved for the    points in $\p D$ and then shown for all other  points, using an argument based on  the Maximum Principle for harmonic functions.  \par
\smallskip
The starting point for Lemma \ref{mainLemma2} is given  by  a preliminary result, which shows that  the constraints  (A) and (D) are satisfied if and only  
if the  vector fields $X_t$ of  an abstract fundamental pair  must have a very special form,  with very  few degrees of freedom. 
With this  the proof boils down  to showing the existence of  abstract pairs $(X_t, J_t)$, in which  the
$X_t$ have the above mentioned special form  and $J_t$ has to satisfy  the  remaining  constraints (B) and (C). Since  (C) is  actually a consequence of (B),  everything reduces to proving the existence of a solution to the differential problem  (B) with $X_t$ in special form. Note that  the explicit expression for a vector field $X_t$ in  special form   involves the complex structure  $J_t$ and the curve $z_t$. This makes the differential equation in (B)  {\it non-linear} in the tensor field $J_t$. 
The proof is  obtained  by first  proving the  existence of solutions  in the  real-analytic category and then getting the   result  by an approximation argument. 
\par
\medskip
An immediate  consequence of the above two lemmas is the fact  that,  for each straight segment $z_t$ in a disk of a Monge-Amp\`ere foliation of sufficiently high regularity, 
there is a one-parameter family of  Monge-Amp\`ere  exhaustions $\t^{(t)}: \overline D {\to} [0,1]$, centered at the points $z_t$. Thus,  the  next  theorem follows. 
\begin{theo}[Propagation of regularity] \label{propofregularity}
Let $\overline D \subset \bC^n$ be a closed strongly pseudoconvex domain of circular type, with  a Monge-Amp\`ere exhaustion $\t: \overline D \longrightarrow [0,1]$, with center $z_o$ and of class $\cC^{k,\a}$ on $\overline D \setminus \{z_o\}$ for some  $k \geq 4$ and $\a > 0$.
Then  there is an open neighborhood $D' \subset D$ of $z_o$ such that any other point $z$ of  $D'$ is center of a Monge-Amp\`ere exhaustion 
$\t^{(z)}{:}\overline D \to [0,1]$
of class   $\cC^{k-2,\a}$ on $\overline D \setminus \{z\}$.
 The dependence of   $\t^{(z)}$ on $z$ is   $\cC^{k-2,\a}$.
\end{theo}
Combining this    with the previously described properties  of Kobayashi metric of domains of  domains of circular type,  we immediately get the following 
\begin{cor} \label{maincor} If $\overline D$ admits a Monge-Amp\`ere  exhaustion, which is $\cC^\infty$ off the center,  its Kobayashi pseudo-metric $\k: TD {\to} [0, + \infty)$ is a  smooth Finsler metric on an appropriate open subset $D' \subset D$. 
\end{cor}
Note that the proof provides also an explicit description of the points of $\p D'$, thus a method to determine when $D' = D$.
From this we see that  the class of domains in $\bC^n$, for which the Kobayashi pseudo-metric is actually a  Finsler metric not only  contains all   smoothly bounded, strictly linearly convex domains, as   Lempert proved, but  it is much larger than that.  It includes for instance  all  smoothly bounded, strongly pseudoconvex circular domains satisfying appropriate conditions and, further, all  sufficiently small deformations  of   such domains (\cite{BD, PS2b}).  This  is indeed a consequence of the property  that
any smooth deformation $\overline D''$  of  a given  circular domain $\overline D$  is 
 biholomorphic to an (abstract)  manifold with boundary,  given   by  equipping  $\overline D$  with  an appropriate deformed  complex structure  $J \neq J_o$.   If the  deformation is sufficiently small,  stability properties  of the equations for stationary disks imply  that $(\overline D, J)$ (hence,  also  $\overline D''$)  admits a   special foliation, made of $J$-stationary disks passing  through  $x_o = 0$  (see e.g. \cite{PS2b}, Prop. 3.4). Then, the regularity of the data  and the stability property of the condition which give $D' = D$ imply that  these disks form  the  Monge-Amp\`ere foliation of an exhaustion $\t: \overline D'' \to [0,1]$, which makes $\overline D''$  a domain of circular type to which Corollary \ref{maincor} applies.  \par
\medskip 
We conclude stressing that all  proofs in \cite{PS3} are actually given  in  the category of complex manifolds with boundary and   are  valid  also for abstract strongly pseudoconvex  manifolds, regardless of their embeddability in $\bC^n$. \par 
\medskip
\section{Concluding remarks}
The above described results can be taken as starting points for various lines of further investigation. Some of them can be shortly described as follows.\par
\smallskip
\noindent(1)  By Theorem \ref{propofregularity}, the existence of a single  Monge-Amp\`ere exhaustion $\t_o$ of class $\cC^{k, \a}$,  $k\geq 4$, $\a > 2$,  off the center implies the existence of an  infinity of  other  Monge-Amp\`ere exhaustions of lower regularity, smaller by two orders.  Such loss of regularity is due  to  a  technical tool used in the proof, namely  the use of  the so-called {\it normalization maps}.  Other than this, we do not see any intuitive reason for  such loss of regularity and we expect that the main results can be  refined on this aspect.  Any   such improvement   would be quite valuable.\par
  Note  also   that if $\t_o$ is a  Monge-Ampere exhaustion, then the   logarithm $u_o = \log \t_o$  is a pluricomplex Green function, with  pole at the  center of $\t_o$.  Hence,  improvements of our results in the  described direction might   give useful information  on pluricomplex Green functions and enrich    the theory   of such  functions that has been so far developed (for the known regularity properties of pluricomplex Green functions,  see e.g. \cite{De, Gu, Bl}).\par
\medskip
\noindent(2)  Corollary \ref{maincor} implies that  if $\overline D$  is a closed domain of circular type  with a smooth Monge-Amp\`ere  exhaustion and satisfying appropriate conditions, then $\overline D$ is completely determined by the Finsler invariants of its Kobayashi metric, namely by  the Finsler curvature and all Finsler covariant derivatives up to an appropriate  order (\cite{Sp}).  On the other hand,  the same  domain is  equipped also with another important sequence of   invariants, the  {\it Bland and Duchamp invariants},  which are tensor fields  on the blow up of $\overline D$ at the  center of a  Monge-Amp\`ere exhaustion, which describe   how  the CR structures of   the level sets $\{\t = c\}$ evolve   when  $c$  varies between $0$ and $1$ (\cite{BD,PS}).  Also these  invariants completely determine $\overline D$ up to biholomorphisms. \par
We feel that it is possible to determine explicit relations between the Finsler  and the Bland and Duchamp invariants. Finding  such relations would very likely lead to a  deep  insight on the intrinsic properties of strictly linearly convex domains and, more generally, of  all domains of circular type. \par
\medskip
\noindent(3) From Lempert theory,  we know  that on any closed,  smoothly bounded,  strictly  linearly convex domain $\overline D$, the disks of the Monge-Amp\`ere foliations  are {\it complex geodesics} for the Kobayashi metric of the domain. For  other types of closed  domains of circular type  with  smooth Monge-Amp\`ere exhaustions,  the disks of the Monge-Amp\`ere foliations are surely  {\it extremal disks}  for the Kobayashi metric at the centers but  there is no manifest reason for them to be  complex geodesics. Since the property of being a complex geodesic can be nicely described  in terms of   Finsler covariant derivatives (see e.g. \cite{AP}), writing down the explicit relations between the Bland and Duchamp invariants and the  Finsler invariants might  be  helpful to characterize  the convexifiable domains (i.e. those that are biholomorphic to some strictly linearly convex domain) in terms of their Bland and Duchamp invariants. Combining this with the so far   known  techniques for constructing domains with prescribed Bland and Duchamp invariants (\cite{BD, PS2}), all this would pave the way 
towards a  useful characterization of  convexifiable domains of $\bC^n$. 
\par

\bigskip
\bigskip
\font\smallsmc = cmcsc8
\font\smalltt = cmtt8
\font\smallit = cmti8
\hbox{\parindent=0pt\parskip=0pt
\vbox{\baselineskip 9.5 pt \hsize=3.1truein
\obeylines
{\smallsmc 
Giorgio Patrizio
Dip. Matematica e Informatica ``U. Dini''
Universit\`a di Firenze
\& Istituto Nazionale di
Alta Matematica
``Francesco Severi''
}\medskip
{\smallit E-mail}\/: {\smalltt patrizio@math.unifi.it
}
}
\hskip 0.0truecm
\vbox{\baselineskip 9.5 pt \hsize=3.7truein
\obeylines
{\smallsmc
Andrea Spiro
Scuola di Scienze e Tecnologie
Universit\`a di Camerino
Via Madonna delle Carceri
I-62032 Camerino (Macerata)
ITALY
}\medskip
{\smallit E-mail}\/: {\smalltt andrea.spiro@unicam.it}
}
}

\end{document}